\documentclass[11pt,a4paper]{article}
\usepackage[brazil]{babel}
\usepackage[latin1]{inputenc}
\usepackage[T1]{fontenc}

\usepackage{graphicx}
\usepackage{amsmath}
\usepackage{dsfont}
\usepackage{amsthm}

\newcommand{\N}{\ensuremath{\mathds{N}}}

\newcommand{\R}{\ensuremath{\mathds{R}}}

\newcommand{\p}{\ensuremath{\mathds{P}}}
\newcommand{\E}{\ensuremath{\mathds{E}}}

\newtheorem{teo}{Teorema}
\newtheorem{Def}[teo]{Definição}
\newtheorem{prop}[teo]{Proposição}

\begin{document}

\begin{center}
\Large{Álbum de figurinhas da Copa do Mundo: uma abordagem via Cadeias de Markov}\

\vspace{1 cm}

\large{Leandro Morgado}

\large{IMECC, Universidade Estadual de Campinas}

\vspace{0.6 cm}

\normalsize{16 de maio de 2014}

\end{center}

\vspace{0.6 cm}

\section{Considerações iniciais}

Com a realização da Copa do Mundo de futebol no Brasil, uma antiga brincadeira tem virado febre em crianças e (especialmente!) adultos em diversas regiões do país. Trata-se de colecionar as figurinhas do álbum oficial da Copa.

Aqueles que apreciam essa brincadeira sabem que a sensação de abrir o próximo envelope e vibrar com uma figurinha inédita na coleção é indescritível. Mas, o fato é que, a partir de certo ponto, esse prazer acaba virando frustração.

Isso porque, a medida que o álbum vai se aproximando do final, as figurinhas repetidas se acumulam de uma forma que fica muito mais cômodo recorrer às trocas com os colegas para conseguir o álbum completo.

Apesar de normalmente seguir o caminho mais fácil (recorrendo às trocas), sempre quis saber, desde criança, quais as chances de completar um álbum apenas comprando um certo número de figurinhas.

E essa curiosidade de infância pode ser respondida conhecendo-se um pouco de Cadeias de Markov e Teoria de Probabilidade.

\section{Modelando o problema}

Para discutir matematicamente o problema, vamos considerar um álbum com $n$ figurinhas. Considere também que um colecionador pretende completar o álbum do jeito mais difícil, ou seja, vai comprando figurinhas de uma em uma, até conseguir o seu objetivo.

\begin{figure}[h!]
\begin{center}
 \includegraphics[scale=0.4]{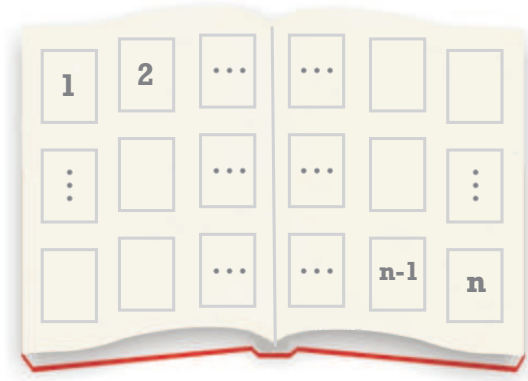}
 \caption{Álbum contendo $n$ figurinhas}
\end{center}
\end{figure}

Finalmente, vamos assumir também que a coleção é ``honesta'', ou seja, que todas as figurinhas tem a mesma probabilidade de serem compradas.

Este problema pode ser modelado por meio de uma sequência de variáveis aleatórias. Vamos denotá-las por $\{X_t: t \geq 0\}$, em que $X_t$ indica o número de figurinhas que o colecionador conseguiu colar em seu álbum após comprar $t$ figurinhas.

Como o total de figurinhas do álbum é $n$, essas variáveis aleatórias podem assumir valores no conjunto $\Omega = \{0, 1, 2, \ldots, n\}$, denominado espaço de estados.

É importante perceber também que o valor assumido por $X_t$ não depende de todo o passado, ou seja, se nas primeiras figurinhas ele deu sorte ou não. Nesse sentido, basta saber quantas figurinhas ele colou no seu álbum após a compra $t-1$, e podemos calcular a probabilidade para a $t$-ésima compra correspondente.

Quando isto acontece, estamos diante de uma cadeia de Markov, assim denominada em homenagem ao matemático Andrei Markov (1856-1922), que obteve resultados importantes sobre fenômenos dessa natureza.

Assim, nesse contexto, para modelar o problema do colecionador de figurinhas, é interessante apresentar a definição a seguir:

\bigskip

\begin{Def}
Seja $\Omega$ um conjunto finito. Uma cadeia de Markov em $\Omega$ é uma sequência de variáveis aleatórias $\{X_t: t \geq 0\}$ com valores em $\Omega$, tal que, se \break$\p \ \big(X_0 = x_0, X_1 = x_1, \ldots, X_t = x\big) \neq 0$, é válida a seguinte propriedade:
$$\p \ \big[X_{t+1} = y \ | \ X_0 = x_0, X_1 = x_1, \ldots, X_t = x\big] = \p \ \big[X_{t+1} = y \ | \ X_t = x\big],$$
para todo $t \in \N$ e para todos $x,y \in \Omega$. Chamaremos o conjunto $\Omega$ de espaço de estados.
\end{Def}

\bigskip

Esta propriedade significa, em termos gerais, que o estado atual do processo depende apenas do estado imediatamente anterior, e não de todo o passado.

Voltando ao problema do Colecionador, podemos representar a transição entre as variáveis aleatórias pelo grafo a seguir:

\vspace{0.5cm}

\begin{figure}[h!]
\begin{center}
 \includegraphics[scale=0.4]{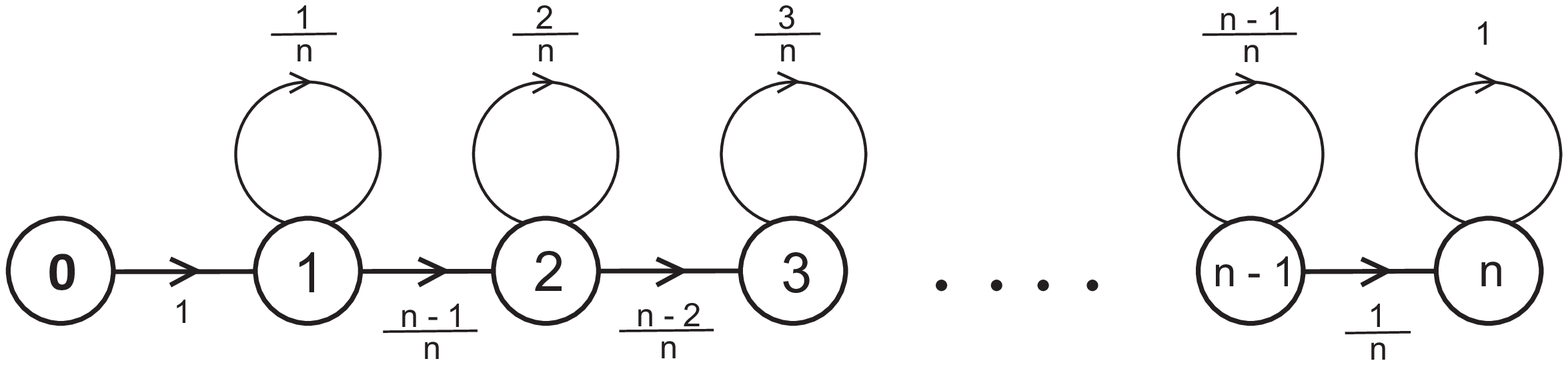}
 \caption{Transição da Cadeia de Markov associada ao problema}
\end{center}
\end{figure}

Uma outra possibilidade para entender o processo é considerar uma matriz quadrada, de ordem $(n+1)$ (número de elementos do espaço de estados), que contém todas as informações sobre a transição do processo.

Nesse sentido, denotando por $A$ essa matriz, cada elemento $a_{ij}$ representa a probabilidade de transição do estado $i$ para o estado $j$.

No nosso problema, temos que a matriz de transição $A$ é dada por:

\vspace{1cm}

\hspace{2cm}$\left(
  \begin{array}{ccccccc}
    0 & 1 & 0 & 0 & \ldots & 0 & 0 \\
    0 & \frac{1}{n} & \frac{n-1}{n} & 0 & \ldots & 0 & 0 \\
    0 & 0 & \frac{2}{n} & \frac{n-2}{n} & \ldots & 0 & 0 \\
    & \vdots &  & &  & & \\
    & \vdots &  & &  & & \\
    0 & 0 & 0 & 0 & \ldots  & \frac{n-1}{n} & \frac{1}{n} \\
    0 & 0 & 0 & 0 & \ldots & 0 & 1 \\
  \end{array}
\right)$

\vspace{0.5cm}

Note que, como não há chances do colecionador perder as figurinhas que já colou em seu álbum, trata-se de uma matriz triangular superior.

\section{Estimativas iniciais do número de figurinhas}

Agora, vamos definir uma variável aleatória que trata do objetivo do colecionador, ou seja, relacionada ao instante em que o este consegue completar o álbum. Nesse sentido:
$$\tau_{fig} = \min \ \{t\geq 0 \ | \ X_t = n\}.$$

Queremos estimar a probabilidade do colecionador ter completado o álbum até um certo tempo $t$. Nesse contexto, o resultado a seguir é fundamental.

\bigskip

\begin{prop}
\label{colecionador}
Na Cadeia de Markov definida acima, tomando $c>0$, temos:
$$\p \ ( \tau_{fig} > \lceil n \ln{n} + cn \rceil ) \leq e^{-c},$$
em que para $a \in \R$, $\lceil a \rceil$ representa o menor inteiro maior que $a$.

\begin{proof}
Note inicialmente que, como assumimos que todas as figurinhas tem a mesma probabilidade de serem compradas, a probabilidade de comprar a $k$-ésima figurinha em um determinado instante é $\frac {1}{n}$. Obviamente, a probabilidade de não comprar a $k$-ésima figurinha nesse mesmo instante é dada por $(1 - \frac {1}{n})$.

Nesse sentido, definimos o evento:
$$A_k = \{\text{k-ésima figurinha não estar entre as  } \lceil n \log{n} + cn \rceil \text{  primeiras}\}.$$

Temos portanto que:
$$P (A_k) = \left( 1 - \frac{1}{n} \right)^{\lceil n \log{n} + cn \rceil}.$$

Note ainda que a probabilidade do álbum não estar completo até o instante $\lceil n \log{n} + cn \rceil$ é dado pela união dos eventos $A_k$, com $k$ variando de 1 a $n$. Dessa forma, podemos escrever:
\begin{equation}\begin{array}{lll}
\p \ ( \tau_{fig} > \lceil n \log{n} + cn \rceil ) &=& \p \left( \displaystyle \bigcup_{k=1}^n A_k \right) \leq \displaystyle \sum_{k=1}^n \p(A_k) = \\
&=& n \cdot \left( 1 - \frac{1}{n} \right)^{\lceil n \log{n} + cn \rceil} = \\
&=& n \cdot \left( 1 - \frac{1}{n} \right)^{n \cdot \frac{\lceil n \log{n} + cn \rceil}{n}} \leq \\
&\leq_{(*)}& n \cdot e^{\frac{-\lceil n \log{n} + cn \rceil}{n}} \leq \\
&\leq& n \cdot e^{\frac{- n \log{n} - cn}{n}} = \\
&=& \frac {n}{n \cdot e^c} = e^{-c},
\end{array}
\end{equation}
em que a desigualdade $\leq_{(*)}$ ocorre em função da convergência monótona $\left(1 - \frac{1}{n}\right)^n \nearrow e^{-1}$ quando $n \to \infty$.
\end{proof}
\end{prop}

O resultado anterior significa, em termos gerais, que se tomarmos o tempo na forma $\lceil n \log{n} + cn \rceil$, a probabilidade do colecionador não ter completado o álbum de figurinhas decai exponencialmente com o aumento de $c$.

\section{Distribuição Geométrica e trocas de figurinhas}

Intuitivamente, é fácil perceber que a medida que estamos próximos de completar o álbum, vai ficando mais e mais difícil conseguir uma figurinha que ainda não possuímos.

Esse fato pode ser confirmado pela Cadeia de Markov descrita anteriormente. Por exemplo, para a probabilidade de obter a última figurinha para completar a coleção é $\frac{1}{n}$ em cada tentativa. E se o número $n$ (total de figurinhas do álbum) for grande, esse valor é de fato muito pequeno.

Nesse contexto, podemos nos perguntar o seguinte:

\begin{enumerate}
\item Em média, de quantas tentativas precisamos para obter esta última figurinha?
\item E se faltarem exatamente $s$ figurinhas, de quantas tentativas (também em média) precisamos para colar mais uma no álbum?
\end{enumerate}

Essas perguntas podem ser respondidas facilmente se conhecermos uma distribuição de probabilidade denominada Geométrica.

Intuitivamente, dizemos que uma variável aleatória $Y$ tem distribuição Geométrica se $Y$ conta o número de tentativas necessárias para se conseguir o primeiro sucesso em eventos independentes e identicamente distribuídos.

Em relação a primeira pergunta acima, esses eventos podem ser descritos como ``comprar a próxima figurinha''. Se esta corresponde à última que falta para completar o álbum, temos um sucesso. Caso contrário, temos um fracasso. Ora, a probabilidade de sucesso em cada evento é $p = \frac{1}{n}$.

Assim, via Distribuição Geométrica, podemos calcular a probabilidade de conseguirmos a última figurinha em exatamente $k$ tentativas. Isso ocorre quando obtemos $(k-1)$ fracassos nas primeiras tentativas e sucesso na $k$-ésima tentativa. Como os eventos são independentes, temos:
$$\p \big[Y = k\big] = (1-p)^{k-1} p.$$

E para obter, em média, o número de figurinhas que precisamos comprar para obter a última, basta calcular a esperança (valor médio) dessa variável aleatória dessa variável aleatória.

Nesse sentido, tomando  $q = p-1$, temos que:

\begin{equation}\begin{array}{lll}
\E [Y] &=& \displaystyle \sum_{k=1}^{\infty} k q^{k-1} p = \\
&=& p \displaystyle \sum_{k=1}^{\infty} k q^{k-1} = \\
&=& p \cdot \frac{d}{dq} \left(\displaystyle \sum_{k=0}^{\infty} q^n \right) = \\
&=& p \cdot \frac{d}{dq} \left(\frac{1}{1-q} \right) = \frac{p}{(1 - q)^2} = \frac{1}{p}.
\end{array}
\end{equation}

De fato, os resultados acima valem para qualquer variável aleatória com distribuição Geométrica, onde $p$ é a probabilidade de sucesso do evento.

Assim, respondendo à primeira pergunta Dessa forma, o número de figurinhas que devem ser compradas para conseguir a última é, em média, $\frac{1}{p} = \frac{1}{1/n} = n$, o que pode ser muito frustrante se o álbum tiver muitas figurinhas.

No mesmo sentido, respondendo à segunda pergunta do início da seção, temos que se faltarem $k$ figurinhas, um sucesso (conseguir mais uma para colar no álbum) tem probabilidade $\frac{k}{n}$. E, portanto, o número médio de figurinhas necessárias para conseguir mais uma é $\frac{n}{k}$.

Isso explica matematicamente porque razão as trocas de figurinhas tem se tornado tão populares (e também econômicas). A medida que o álbum vai ficando mais e mais completo, a relação $\frac{n}{k}$ vai crescendo, ou seja, seria preciso comprar muito mais figurinhas para conseguir a próxima.

\section{Estimativas sobre o álbum oficial da Copa}

Nesta seção, vamos aplicar os resultados anteriores para obter estimativas concretas sobre o álbum oficial da Copa do Mundo.

Para este caso específico, o número total de figurinhas é  $n = 649$ (times e estádios numerados de 1 a 639, além das figurinhas 00, W1, J1, J2, J3, J4, L1, L2, L3 e L4).

Pela proposição \ref{colecionador}, obtemos o seguinte resultadö:
$$\p \ ( \tau_{fig} > \lceil n \ln{n} + cn \rceil ) \leq e^{-c}.$$

Nesse sentido, se um colecionador quisesse saber o número de figurinhas necessárias para completar o álbum com certeza, a resposta seria ``não existe esse número''.

Obviamente, por se tratar de um processo aleatório, somente podemos estimar o número de figurinhas necessário para que ele complete o álbum com uma determinada probabilidade, e é nesse sentido que vamos usar o resultado acima.

Assim, para fixar as ideias, vamos estimar essa probabilidade em 90\%. Agora, vamos estimar uma cota para o número de figurinhas que deveriam ser compradas, a fim de que a probabilidade de não completar o álbum seja inferior a 10\%.

Para tanto, tome $c$ tal que $e^{-c} \approx 0,1$. Temos então $c = - \ln 0,1 \approx 2,3$. Fazendo as contas, temos que $\lceil 649 \cdot \ln{649} + 2,3 \cdot 649 \rceil = 5696$, e portanto:
$$\p \ ( \tau_{fig} > 5696 ) \leq 0,1.$$

Ou seja, se o colecionador comprar \textbf{5696 figurinhas (!!!)}, podemos garantir que a chance de completar o álbum da Copa é superior a 90\%. Como o preço de cada figurinha é R\$ 0,20, estamos falando de um investimento total de \textbf{R\$ 1.139,20}.

Finalmente, para obter a última figurinha do jeito mais difícil, vimos, na seção referente à distribuição Geométrica, que precisamos, em média, de $n = 649$ tentativas. Isso corresponde (em média) a R\$ 129,80.

Vendo essas estimativas, percebemos que, de fato, é bem melhor recorrer às trocas, não é mesmo?

\section{Referências}

\begin{itemize}
\item FELLER, W. {\it An Introduction to Probability Theory and its Applications}. vol. 1. New York, 1993.

\item JAMES, Barry. {\it Probabilidade: um curso em nível intermediário}. Rio de Janeiro: Instituto de Matemática Pura e Apicada, 2010.

\item MORGADO, Leandro. {\it Tempos de Mistura em Cadeias de Markov}. Trabalho de Conclusão de Curso. Universidade Federal de Santa Catarina. Florianópolis, 2011.

\item PERES, D.L.Y. {\it Markov Chains and Mixing Times}. American Mathematical Society, 2008.
\end{itemize}

\end{document}